\documentclass[review]{elsarticle}

\usepackage{lineno,hyperref}
\modulolinenumbers[5]

\usepackage{amsmath}
\usepackage{amsfonts}
\usepackage{amssymb}
\usepackage{amsthm}
\usepackage{mathtools}
\usepackage{algorithm}

\usepackage{algpseudocode}
\usepackage{graphicx}
\usepackage{caption}
\usepackage{subcaption}
\usepackage{multicol}
\usepackage{multirow}
\usepackage{hhline}
\usepackage{float}
\usepackage{color}
\usepackage{xcolor}
\usepackage{chngcntr}
\usepackage{tikz}
\usepackage{array}
\usepackage{comment}
\usepackage{geometry}
\usepackage{bm}
\usepackage{epsfig}

\usepackage{latexsym}
\usepackage{enumerate}
\usepackage{hyperref}
\usepackage{cleveref}

\usepackage[normalem]{ulem}
\usepackage{bbm}
\usepackage{commath}

\graphicspath{{./figs/}}

\newcommand{\ds}{\displaystyle}

\def\wb{\overline}
\def\wh{\widehat}

\def\qan{{\quad\hbox{and}\quad}}

\newcommand{\0}{{\mathbf{0}}}

\newcommand{\m}{{\mathbf{m}}}

\newcommand{\bu}{{\mathbf{u}}}
\newcommand{\bv}{{\mathbf{v}}}
\newcommand{\by}{{\mathbf{y}}}

\newcommand{\bbR}{\mathbb{R}}

\newcommand{\bbZ}{\mathbb{Z}}

\newcommand{\bbeta}{{\boldsymbol\eta}}
\newcommand{\bzeta}{{\boldsymbol\zeta}}
\newcommand{\bPsi}{{\boldsymbol\Psi}}
\newcommand{\btheta}{{\boldsymbol\theta}}

\newcommand{\cN}{\mathcal{N}}

\newcommand{\cF}{\mathcal{F}}
\newcommand{\cH}{\mathcal{H}}

\numberwithin{equation}{section}
\allowdisplaybreaks

\newtheorem{remark}{Remark}[section]

\usepackage[normalem]{ulem}

\begin{document}

\begin{frontmatter}

\title{Non-intrusive Structural-Preserving Sequential Data Assimilation}


\author[dartmouth]{Lizuo Liu\fnref{myfootnote}\corref{mycorrespondingauthor}}
\cortext[mycorrespondingauthor]{Corresponding author}
\ead{lizuo.liu@dartmouth.edu}

\author[umbc]{Tongtong Li\fnref{myfootnote}\fnref{myfootnote2}}
\ead{tongtong.li@umbc.edu}

\author[dartmouth]{Anne Gelb\fnref{myfootnote}}
\ead{annegelb@math.dartmouth.edu}

\address[umbc]{Department of Mathematics and Statistics, University of Maryland, Baltimore County, Baltimore, MD 21250, USA}

\address[dartmouth]{Department of Mathematics, Dartmouth College, Hanover, NH 03755, USA}

\fntext[myfootnote]{This work is partially supported by the NSF grant DMS \#1912685, DOE ASCR \#DE-ACO5-000R22725, and DOD ONR MURI grant \#N00014-20-1-2595.}

\begin{abstract}

Data assimilation (DA) methods combine model predictions with observational data to improve state estimation in dynamical systems, inspiring their increasingly prominent role in geophysical and climate applications. Classical DA methods assume that the governing equations modeling the dynamics are known, which is unlikely for most real world applications. Machine learning (ML) provides a flexible alternative by learning surrogate models directly from data, but standard ML methods struggle in noisy and data-scarce environments, where meaningful extrapolation requires incorporating physical constraints. Recent advances in structure-preserving ML architectures, such as the development of the entropy-stable conservative flux form network (ESCFN), highlight the critical role of physical structure in improving learning stability and accuracy for unknown systems of conservation laws. Structural information has also been shown to improve DA performance. Gradient-based measures of spatial variability, in particular, can help refine ensemble updates in discontinuous systems. Motivated by both of these recent innovations,  this investigation  proposes a new non-intrusive, structure-preserving sequential data assimilation (NSSDA) framework that leverages structure at both the forecast and analysis stages. We use the ESCFN to construct a surrogate model to preserve physical laws during forecasting, and a structurally informed ensemble transform Kalman filter (SETKF) to embed local statistical structure into the assimilation step, assuming partial structural information such as the transformation between physical and conserved variables. Our method operates in a highly constrained environment, using only a single noisy trajectory for both training and assimilation. Numerical experiments where the unknown dynamics correspond respectively to the shallow water and Euler equations demonstrate significantly improved predictive accuracy.

\end{abstract}


\end{frontmatter}


\section{Introduction}\label{sec:introduction}

Data assimilation (DA) provides a statistically accurate and reliable computational framework for integrating numerical model predictions with observational data to obtain system state estimates. It has been broadly used in applications such as climate modeling \cite{Hannart16}, weather forecasting \cite{Kalnay02}, and sea ice dynamics \cite{Cheng23}. Within the context of DA, {\em filtering} describes a methodology that sequentially updates the system state estimates as new observations become available over time, while {\em smoothing} retrospectively improves past state estimates by incorporating both earlier and future measurements. This investigation is concerned with discrete-time filtering from a probabilistic perspective, that is, for which the conditional distribution of the state variable is sequentially updated based on all data accumulated up to the current time in a two-step procedure.  In the forecast (or prediction) step, the system’s dynamical model is used to propagate the current distribution forward, producing a prior distribution at the next time step. In the subsequent analysis step, Bayes’ theorem is applied to incorporate newly acquired observations (interpreted as a likelihood) into the prior, resulting in a posterior distribution. This posterior reflects all information available so far and serves as the updated estimate of the system’s state. While smoothing offers valuable complementary capabilities, it falls outside the scope of the present study.


Governing models in realistic applications are typically highly nonlinear and therefore violate the assumptions underlying  classical Kalman filters  \cite{Kalman60}, which are optimal only for linear dynamics and Gaussian noise \cite{Cohn97}.  Ensemble-based Kalman filters \cite{Evensen94} provide a robust and computationally efficient alternative by approximating the state distribution via a set of ensemble members and using empirical covariances to perform updates, thereby enabling practical solutions for high-dimensional systems. However, the use of low-rank estimates can fail to capture important physical and structural features, particularly in geophysical applications where sharp gradients or discontinuities (e.g., fronts, shocks, or ridges) are common. This limitation was successfully addressed in \cite{TLda}, where  the structurally-informed ensemble transform Kalman filter (SETKF) framework was introduced.  Specifically, SETKF replaces the standard covariance (global) weighting matrix which is used to balance the prior and observational information with a structurally-informed weighting matrix constructed from the second moment of the state variable's discrete gradient. SETKF is  therefore able to capture spatial (local) variability and emphasize regions of structural change, consequently improving the accuracy and stability of the assimilation updates.

A remaining significant challenge is that standard ensemble filtering techniques typically require full knowledge of the governing equations and system parameters. This assumption often fails in geophysical applications, where models may be incomplete, inaccurate, or entirely unavailable. Furthermore, sometimes only a qualitative understanding of the underlying physics may be accessible, that is, no precise mathematical model can be formulated as a system of ordinary or partial differential equations. For example, while it may be known that sea ice dynamics involve advection and satisfy certain conservation laws, the exact form of these laws or the relevant transport velocities may not be readily apparent. In such cases machine learning (ML) methods serve as attractive purely data-driven alternatives, capable of learning effective dynamical models by minimizing prediction errors over one or several time steps based on noisy observations. However, in the absence of explicit physical constraints, ML models trained on noise-contaminated data are prone to instability and can yield unphysical long-term predictions. Physics-informed neural networks (PINNs) \cite{Pinns2019} introduce soft physical constraints into the training objective, and while improving short-term performance, the methodology often remains insufficient for stabilizing predictions over long-term horizons. More recently ML architectures have been designed with explicit structural priors tailored to the physics of the problem. Notable examples include the Conservative Flux Form Network (CFN) \cite{chencfn}, the Entropy-Stable CFN (ESCFN) \cite{Liu2024}, and the Neural Entropy-Stable CFN (NESCFN) \cite{Liu2025NESCFN},
which are specifically constructed to respect the conserved form and entropy conditions of hyperbolic conservation laws. These structure-preserving architectures have demonstrated substantial improvements in learning stability and predictive accuracy for systems with complex dynamics.

The aforementioned developments in structure-preserving machine learning make it natural to consider combining ML with DA, particularly in model-free settings where the governing equations are not available and ML must serve as a surrogate. To be clear, there is already growing interest in integrating ML and DA, although current approaches have different target objectives, and also assume at least partial knowledge of the underlying physical model. Such hybrid approaches may be broadly categorized into three groups: The first involves applying ML models as corrective mechanisms to augment imperfect physical models, thereby enhancing the accuracy of the posterior estimates \cite{mlerrorcorection, deepDataAssimilation}. The second utilizes ML methods  (particularly autoencoders) as model reduction tools so that high-dimensional state variables may be projected into a more efficient low-dimensional latent space.  In addition to reducing computational complexity, methods in this category are also designed to enable high-resolution inpainting of the full state space from latent representations \cite{datalearning}. Finally, algorithms in the third category draw inspiration from generative modeling, especially diffusion-based methods, where the posterior distribution at each assimilation step is approximated via a learned score function conditioned on the available observations \cite{scorefilter, DiffDA}. 

We emphasize that methods in each of these categories are generally built on partially known dynamics and do not explicitly leverage structural information. In contrast, our approach centers on the role of structural information, both {\em physical} (e.g., conservation laws, entropy conditions) and {\em statistical} (e.g., spatial gradient variability), as the primary means of learning and assimilating in the absence of any governing equations. We consider a purely non-intrusive and data-driven setting in which only a single trajectory of noisy observations is available, which is again in contrast to prior work (e.g., \cite{Liu2024}) that relies on multiple noise-free trajectories for training. This one trajectory is used both for learning a surrogate model as well as conducting sequential assimilation, reflecting practical scenarios where prior information is limited and uncertainty is high. To compensate for this scarcity of data and absence of models, we exploit structural information in {\em both} components of the assimilation pipeline. On the one hand, structural priors inform the learning of a surrogate model that is not exact but preserves key physical properties. On the other hand, incorporating statistical structure into the assimilation step, which can be viewed as using observational data to correct accumulated forecast error,  enables more accurate and targeted updates. 

To this end we propose a new non-intrusive, structure-preserving sequential data assimilation (NSSDA) framework that combines ESCFN with SETKF. Our method enforces structural consistency throughout the pipeline: the surrogate model learned by ESCFN preserves physical principles through entropy-stable and conservative design, while the SETKF analysis step leverages second-moment information of the discrete gradient of the prior to guide updates in regions of sharp features or discontinuities. This dual incorporation of physical and statistical structure enables the framework to extract predictive information from limited and noisy data, even in the absence of a known governing model. Since ESCFN is trained to predict the dynamics of conserved variables, and SETKF performs analysis updates on the physical variables using observational data, a known mapping from physical to conserved quantities is assumed. In this sense, our method is partially physics informed. We emphasize that this assumption is significantly weaker than assuming knowledge of the full state or dynamical model. For instance, one might measure pressure and know it contributes to a conserved energy variable, even without knowing the exact conservation law or its form.
We demonstrate our new methodology on  (1) the shallow water equations and (2) the Euler equations.  Our results illustrate the effectiveness of structure-aware learning and assimilation in recovering complex dynamics from minimal information, highlighting the broader potential of structure-preserving methods in data-driven scientific inference.

The rest of the paper is organized as follows. In Section~\ref{sec:preliminaries} we review the foundational components of our framework, including the ensemble transform Kalman filter (ETKF), its structurally informed extension SETKF, and the surrogate modeling architectures including Neural ODE and the ESCFN. Section~\ref{sec:method} describes our proposed NSSDA approach, which is specifically designed to enforce structural consistency by coupling ESCFN-based surrogate modeling with SETKF-based data assimilation. Numerical experiments on the shallow water and Euler equations are presented in Section~\ref{sec:numerical} to illustrate the benefits of incorporating structural information in both the forecasting and assimilation stages. We conclude with a discussion of findings and future directions in Section~\ref{sec:conclusion}.

\section{Preliminaries}
\label{sec:preliminaries}

We begin by reviewing the foundational components of our proposed {\em Non-Intrusive Structure-Preserving Sequential Data Assimilation} (NSSDA) framework. NSSDA combines structure-preserving surrogate modeling with structure-informed ensemble-based data assimilation. To set the stage, we first review the key building blocks, specifically the classical ETKF and its structurally informed variant SETKF for assimilation, along with two data-driven forecasting models: (1) the standard Neural ODE \cite{neuralode} and (2) the ESCFN \cite{Liu2024}.

We consider a standard sequential data assimilation setting, in which the discrete-time state variables $\{\bv_j\}$ evolve according to some dynamical process, and noisy observational data $\{\by_j\}$ are collected at discrete time steps. That is, the dynamical model is of the general form
\begin{equation}\label{eq: dynamic}
\bv_{j+1}=\bPsi(\bv_j), \ j = 0, \cdots, J-1,
\end{equation}
where $\bPsi \in C( \bbR^n, \bbR^n)$ denotes a deterministic evolution operator. This operator typically arises from a numerical discretization of an underlying continuous dynamical system governed by some (known) differential equations. For ODEs, this involves applying a time integration scheme directly. For PDEs, the process typically involves first discretizing the spatial domain, for example, via finite difference or finite element methods, which is then followed by temporal integration of the resulting semi-discrete system. The time instances are given by $t_j=j \Delta t$, where $\Delta t = T/J$, and $J \in \bbZ^{+}$ is the final time step for the simulation. The initial state $\bv_0$ is assumed to follow a Gaussian distribution 
\begin{equation}\label{eq: init Gaussian}
  \bv_0 \sim \cN( \m_0, C_0).
\end{equation}

The second part of sequential data assimilation involves the given data or observations $\{\by_j\}$, which are modeled as a linear function of the signal with  additive noise $\{\bbeta_{j}\}$ as
\begin{equation}\label{model: observation}
\by_{j+1}=H\bv_{j+1}+\bbeta_{j+1}, \ j = 0, \cdots, J-1.
\end{equation}
Here $H\in C(\bbR^n, \bbR^m)$ is the observation operator, and the noise $\{\bbeta_j\}$ is an independent and identically distributed (i.i.d.) sequence, independent of $\bv_0$, with
$\bbeta_1 \sim \cN(\0,\Sigma)$,
where $\Sigma$ is a known, symmetric positive definite (SPD) covariance matrix. In the special case where the observation noise $\bbeta_1$ consists of i.i.d.~random variables with variance $\sigma^2$, the covariance matrix simplifies to $\Sigma = \sigma^2 I$, with $I$ being the identity matrix.

In the setting of interest, the governing equations that define $\bPsi$ are assumed to be {\em unavailable or unknown}. Specifically, this work focuses on systems governed by hyperbolic conservation laws, where full knowledge of the underlying model is unavailable, but partial structural information, such as the conversion between physical and conserved variables (e.g., velocity and pressure to momentum and energy), is known. The operator can be approximated by a data-driven surrogate model that is learned from the observational data, such as one constructed using neural networks. The goal is to reconstruct the latent state from noisy observations of physical state variables using a non-intrusive sequential data assimilation framework. A more detailed formulation for the hyperbolic conservation law setting is provided in \Cref{subsec: hyperbolic conservation law}. 


\subsection{Ensemble transform Kalman filter (ETKF)}
\label{subsec: ETKF}
When the system dynamics (modeled by $\bPsi$) are linear and all associated distributions are Gaussian, the classical Kalman filtering method \cite{Kalman60} yields an optimal sequential update of the state probability distribution conditioned on data, as the posterior remains Gaussian and is entirely described by its mean and covariance \cite{Law15, Cohn97}. In most practical applications, however, the system dynamics are typically nonlinear, and the Gaussian assumptions break down.  Various ensemble-based Kalman filter methods have been developed as practical alternatives, including the Ensemble Kalman Filter (EnKF) \cite{Evensen94}, the Ensemble Adjustment Kalman Filter (EAKF) \cite{Anderson01}, and the Ensemble Transform Kalman Filter (ETKF) \cite{Bishop01}. The idea is to sequentially estimate the state distribution by propagating an ensemble of particles through the model and updating them using observations. In what follows we briefly describe the ETKF method as it serves as a prototype for our approach.

The ETKF method is realized in two steps: forecast and analysis. In the forecast step, the particle approximations $\{ \bv^{(k)}_j\}_{k=1}^{K}$ of the filtering distribution are propagated through the transition dynamics to yield the forecast ensembles $\{ \wh{\bv}^{(k)}_{j+1} \}_{k=1}^{K}$ at the next observation time $t_{j+1}$ for each $k$ as
\begin{equation}\label{eq: EnKF predict 1}
\wh{\bv}^{(k)}_{j+1}=\bPsi(\bv^{(k)}_j).
\end{equation}
The forecast distribution is then approximately characterized by its corresponding sample prior mean and covariance matrix
\begin{subequations}
\begin{equation}
\ds \wh{\m}_{j+1}=\frac{1}{K}\sum_{k=1}^{K}\,\wh{\bv}^{(k)}_{j+1}, \label{eq: EnKF prior 1}
\end{equation}
\begin{equation}
\ds \wh{C}_{j+1} =\frac{1}{K-1} \sum_{k=1}^{K}\,(\wh{\bv}^{(k)}_{j+1}-\wh{\m}_{j+1})(\wh{\bv}^{(k)}_{j+1}-\wh{\m}_{j+1})^{T},  \label{eq: EnKF prior 2}   
\end{equation}
\end{subequations}
which respectively serve as low-dimension representations of the prior mean and covariance. The forecast ensembles $\{ \wh{\bv}^{(k)}_{j+1} \}_{k=1}^{K}$ are then updated to obtain posterior ensembles in the analysis step using new observation data $\by_{j+1}$. This is accomplished by first updating the prior mean to the posterior mean via minimization as
\begin{equation}
\m_{j+1} = \underset{\bv}{\text{arg min}} 
\left(\frac{1}{2}\vert \by_{j+1}-H \bv \vert_{\Sigma}^2 + \frac{1}{2} \vert \bv-\wh{\m}_{j+1}^{(k)} \vert_{W_{j+1}}^2\right).
\label{eq:ETRF min}
\end{equation}
In traditional ETKF the weighting matrix $W_{j+1}$  is given by the $\wh{C}_{j+1}$ in \eqref{eq: EnKF prior 2}, while the posterior covariance matrix is obtained through the Kalman update formula
\begin{equation}
C_{j+1}=(I-K_{j+1}H) \wh{C}_{j+1}, \label{eq:cov_update}
\end{equation}
where the Kalman gain $K_{j+1}$ and the innovation covariance $S_{j+1}$ are respectively given by
\begin{subequations}
\begin{equation}
K_{j+1}=\wh{C}_{j+1}H^{T}S_{j+1}^{-1}, \label{eq:gain_update}
\end{equation}
\begin{equation}
S_{j+1} =H \wh{C}_{j+1}H^{T} + \Sigma. \label{eq:S_update}
\end{equation}
\end{subequations}
The posterior ensemble is then realized as 
\begin{equation}
\ds \bv^{(k)}_{j+1} = \m_{j+1}+\bzeta^{(k)}_{j+1},\quad k = 1,\dots,K,
\label{eq:zeta}
\end{equation}
where each $\bzeta^{(k)}_{j+1}$ is constructed using an appropriate transformation operator (details in \cite{Bishop01, TLda}) to ensure a Gaussian distribution $\cN(\0,C_{j+1})$. 

For the practical implementation of ETKF, covariance inflation \cite{Anderson07, Law15} and localization \cite{Evensen06} are commonly used to ensure the effectiveness of the observation and removal of unwanted spurious correlations in the data assimilation process. Localization is implemented by applying a tapering or masking operator to the sample covariance, typically realized through a Hadamard (elementwise) product with a predefined localization matrix. When the observation data is directly obtained and dense, that is, the observation operator $H$ is the identity operator $I$, the localization matrix reduces to the identity, and the localized covariance simplifies to the diagonal of the sample covariance. Multiplicative inflation is then applied to this localized prior covariance to account for potential underestimation of ensemble spread. Altogether, both localization and inflation can be simultaneously realized by setting 
\begin{equation}
W_{j+1} = \alpha^2 \text{ diag}(\wh{C}_{j+1})
\label{W: cov}
\end{equation}
in \eqref{eq:ETRF min}, where $\alpha$ is the inflation parameter.

While $W_{j+1}$ clearly plays a critical role in balancing the model and data information in \eqref{eq:ETRF min}, the sample covariance used in its construction fails to carry structural information in systems containing discontinuous or highly localized features, such as shocks or fronts, which may in turn lead to inaccurate or unphysical state estimates. 

\subsection{Structurally-informed ETKF (SETKF)}
\label{subsec: SETKF}
The structurally-informed ETKF (SETKF) was proposed in \cite{TLda} to overcome the limitations of covariance-based analysis in systems with sharp gradients or discontinuities. In particular, the weighting matrix $W_{j+1}$ in \eqref{eq:ETRF min} is modified to incorporate local structural information derived from the forecast ensembles by utilizing the discrete gradient information. In the dense observation case, this is accomplished  by replacing \eqref{W: cov} with a diagonal matrix that reflects spatial variability, and is given by
\begin{equation}\label{eq: diag gsm}
W = \beta \wh{S},\quad \beta > 0.
\end{equation}
Here $\beta$ is a tuning parameter similar to the inflation parameter $\alpha$ in \eqref{W: cov}, while $\wh{S}$ is a diagonal matrix with nonzero entries
$\wh{S}_{i,i} = \frac{1}{2}(\wh{S}_{i-\frac{1}{2}}+\wh{S}_{i+\frac{1}{2}})$, where 
\begin{equation}
\label{eq:secmomentcenter}
\wh{S}_{i+\frac{1}{2}}=\frac{1}{K} \sum_{k=1}^{K} \left( d\wh{v}^{(k)}_{i+\frac{1}{2}} \right)^2, \qquad d\wh{v}^{(k)}_{i+\frac{1}{2}}=\frac{\wh{v}^{(k)}_{i+1}-\wh{v}^{(k)}_{i}}{\Delta x}.
\end{equation}
That is, the structural information of the system is obtained by computing the second moment of the first order finite differencing gradient of the state variable.
This construction adaptively assigns smaller weights (i.e., greater trust in observations) in regions with large gradient variability, such as near discontinuities, and larger weights (i.e., greater trust in the model forecast) in smoother regions. Through rebalancing the influence of observations and model predictions with respect to variability or discontinuities of the state variable, the result is a more accurate and physically consistent analysis, especially in systems governed by hyperbolic conservation laws. 

While this approach improves robustness and localization in systems with spatial structure,  it has until now  assume that the corresponding forecast ensembles are generated from a meaningful underlying model. In what follows we show how the prediction step may be adapted in the absence of such a model.
\subsection{Surrogate model training}
\label{sec:trainingperiod}

In the context of data assimilation, under the assumption that a temporal sequence of observations $\{\by_j\}_{j=1}^{J}$ is available, the surrogate model can be trained from a subset of this trajectory, e.g., $\{\by_j\}_{j=1}^{L_{train}}$ with $L_{train}\leq J$. Once trained on these observational snapshots, the learned flow map can be used in place of $\bPsi$ during the prediction step of ensemble-based filters to propogate each ensemble member in \eqref{eq: EnKF predict 1}:
\begin{equation}\label{eq: EnKF predict 1 - Neural ODE}
\wh{\bv}^{(k)}_{j+1}=\bPsi^\btheta(\bv^{(k)}_j),
\end{equation}
where $\bPsi^{\btheta}$ denotes the one-step flow map induced by solving the Neural ODE \eqref{eq: neural ODE}, which we write as $\bPsi^\btheta_{NN}$, or ESCFN \eqref{eq:Ndefine}, denoted by $\bPsi_{CFN}^\btheta$, over the interval $[t_j, t_{j+1}]$. This operator is learned exclusively from data, thereby allowing the prediction step to be performed non-intrusively  without explicit knowledge of the underlying governing equations.

\subsection{Neural ODEs}
\label{sub:neuralODE}
We begin by considering the Neural Ordinary Differential Equations (Neural ODEs) \cite{neuralode} as one possible approach to construct a surrogate dynamical model when  $\bPsi$ in \eqref{eq: EnKF predict 1} is either unavailable. or unknown. 
Instead of modeling discrete-time transitions directly, Neural ODEs define the system evolution via an ODE governed by a neural network for conserved variable
\begin{equation}\label{eq: neural ODE}
\frac{d\bu(t)}{dt} = f^\theta(\bu(t), t),
\end{equation}
where $\bu(t)\in \bbR^{n}$ denotes the system state variable at time $t$, and $f^\theta$ is a neural network parameterized by $\theta$. This formulation originally is viewed as a continuous-depth generalization of residual neural networks, where the transformation of the state is defined by solving an initial value problem over a time interval. 

The Neural ODE framework offers several practical advantages. First, it enables adaptive time stepping through black-box ODE solvers and second, it supports memory-efficient training via the adjoint sensitivity method.  These properties make it particularly well suited for modeling irregularly sampled data or systems with inherently continuous dynamics. Furthermore, Neural ODEs have sparked considerable interest and spawned numerous extensions for time-series analysis \cite{neuralodeTimeSeries} and generative modeling \cite{nodeVAE}, highlighting their potential as a flexible and theoretically grounded learning paradigm.

Despite their flexibility, Neural ODEs also face notable limitations. Training can be computationally intensive due to the potentially large number of function evaluations required by ODE solvers. More importantly, these models often lack physical interpretability and may fail to respect structural properties of the underlying system, such as conservation laws, stability, or entropy conditions. This lack of structure can be particularly problematic when modeling systems with discontinuities or sharp gradients, as standard Neural ODEs tend to produce overly smoothed predictions that may violate physical constraints. Such limitations have motivated the development of surrogate architectures that incorporate physical structure directly into their formulation \cite{chencfn, Liu2024, Liu2025NESCFN}. 
Below we review the Entropy-Stable Conservative  Flux Form Network (ESCFN) \cite{Liu2024} which we will then use for the prediction step in our new NSSDA approach.

\subsection{Entropy-Stable Conservative  Flux Form Network (ESCFN)}
\label{sec:ESCFN}

The entropy-stable conservative flux form neural network (ESCFN) developed in \cite{Liu2024} addresses the limitations of unconstrained neural surrogates by incorporating essential physical principles into the learned model. It is designed to learn the system dynamics from data while explicitly preserving core structural properties, including conservation laws and entropy conditions, and when solutions develop shock discontinuities. Specifically, ESCFN mimics the numerical structure of high-resolution finite volume schemes. In so doing, it enables stable predictions aligned with the physical constraints of the governing PDEs, even when trained on noisy data. In contrast to Neural ODEs, which may yield overly smoothed or unphysical solutions, ESCFN ensures stability and sharp feature resolution through conservation and entropy stability. As such, ESCFN plays a central role in the NSSDA framework, enabling accurate and physically consistent predictions from data alone.  The details of ESCFN will be discussed in \Cref{subsec: surrogate model training}.

\section{The non-intrusive structure-preserving sequential data assimilation (NSSDA) method}
\label{sec:method}
The preceding discussion highlights a critical insight: structure must be preserved throughout the entire data assimilation pipeline. Using a structure-informed analysis method with a poorly trained or unstructured forecast model leads to unstable or misleading updates. Conversely, applying a structure-preserving surrogate model without respecting structural variability in the assimilation step can lead to suboptimal data incorporation.
Our new NSSDA approach unifies these ideas of both structure-preserving surrogate modeling and structure-informed assimilation into a single framework. 
The prediction step is handled by a surrogate trained directly from observational data using the structure-preserving ESCFN, and the analysis step leverages the preserved localized structural information to ensure robust updates.
\subsection{Hyperbolic conservation laws and numerical discretization}
\label{subsec: hyperbolic conservation law}
The focus in this investigation is on data assimilation for systems governed by hyperbolic conservation laws. That is, we assume the underlying dynamics satisfy the PDE
\begin{equation}
\dfrac{\partial u}{\partial t}  + \dfrac{\partial f\left( u \right)}{\partial x} = 0, \quad x \in (a,b), \quad t \in (0,T), 
\label{eq: conservation_law}
\end{equation}
with appropriate boundary and initial conditions.
Here $u = u(x,t)$ denotes the conserved state variable and $f(u)$ is the flux function. These equations arise in a wide range of physical systems and are characterized by the possibility of developing discontinuities, such as shock waves, even when starting from smooth initial data.

The explicit form of the PDE  \eqref{eq: conservation_law} is {\em not} assumed to be known. Rather we adopt a non-intrusive approach in which a surrogate model is trained directly from data. As was demonstrated in \cite{chencfn,Liu2024,Liu2025NESCFN}, in order to obtain the correct  speed and location of discontinuous fronts, it is essential to build the surrogate model using a conservative numerical flux form that respects the core properties of the PDE. To that end, we consider the semi-discrete conservation form, which captures the spatial structure of the PDE and serves as the basis for surrogate modeling in our framework.

We begin by discretizing the spatial domain $[a,b]$ with grid points $\{x_i\}_{i = 0}^n$, where
\begin{equation}
\label{eq:spatialgrid}
x_i=a + i \Delta x, \quad \Delta x = \frac{b-a}{n}.
\end{equation} 
We approximate the solution using cell averages
\begin{equation}
\label{eq: cell_average}
\wb{u}_i(t) = \frac{1}{\Delta x}\int_{x_i - \frac{\Delta x}{2}}^{x_i + \frac{\Delta x}{2}} u(x,t)\,dx, \quad i=1, \cdots, n-1, 
\end{equation}
with appropriate boundary conditions applied at $i=0$ and $i=n$. The evolution of the cell averages is governed by the semi-discrete conservation law 
\begin{equation}
\frac{d}{dt} \wb{u}_i + \frac{1}{\Delta x} \left(  f_{i+1/2} - f_{i-1/2} \right) = 0,
\label{eq: semi-discrete conservative scheme}
\end{equation}
where $f_{i\pm 1/2}$ are numerical fluxes that approximate the true flux across cell interfaces $x=x_{i\pm 1/2}$. This conservation form is crucial for maintaining consistency with the underlying PDE and capturing correct shock speeds \cite{Hesthaven,LeVeque92,LeVeque02}.

The semi-discrete system \eqref{eq: semi-discrete conservative scheme} is integrated using explicit time-stepping methods suitable for hyperbolic systems, such as total variation diminishing (TVD) Runge–Kutta methods \cite{Liu94}, strong-stability-preserving Runge-Kutta, which are designed to capture discontinuities while preserving numerical stability \cite{SSPscheme}.

The spatial resolution $\Delta x$ is chosen to ensure that the solution is well-resolved and to guarantee stability. In practice, this resolution will depend on computational feasibility, and will play a key role in the performance and stability of both the surrogate model training and the data assimilation process. The uniform time step $\Delta t = \#CFL \Delta x$ is correspondingly chosen to satisfy the Courant–Friedrichs–Lewy (CFL) stability condition \cite{Courant1928, Courant1967}.
We will also assume that we have available $J$ observations in the temporal interval $[0,T]$, with $J = int(\frac{T}{\Delta t})$.\footnote{For simplicity, we then redefine $T = J\Delta t$ so that we do not have to use any fractional time steps.}


\subsection{Problem setup and data generation}

With the conservation property of underlying dynamics established, we now discuss the data setting and underlying assumptions used to formulate the surrogate learning and data assimilation framework. 

In practical settings the state variable $u$ in \eqref{eq: conservation_law} represents a conserved quantity that may not be directly measurable.  Observations  typically consist of more ``physical'' quantities,  such as density, velocity, or pressure, which relate to the conserved variables through some arithmetic operations, such as multiplications and/or linear combinations. We define such an observable quantity as $v$, and assume the  mapping between $u$ and $v$, denoted by $\Phi : u \to v$, is a known and invertible transformation. We further assume that the observable data consist of a set of discrete temporal snapshots of the physical variables $v$. We mimic the realistic observational constraints in which only limited measurements along one solution path are accessible by generating the snapshots from the numerical solution of a prototype hyperbolic conservation law over an initial time window. 


Importantly, we stress that the flux functions \(f(u)\) in the system corresponding to each conserved quantity are  {\em not} known apriori. Our goal is to approximate the dynamics of the conserved variables $u$ of the hyperbolic system, learn a surrogate that captures its temporal evolution, and ultimately use this surrogate within an ensemble-based data assimilation framework to predict the future state of the physical variables $v$ beyond the initial training window. For ease of presentation we describe our new methodology for the one-dimensional scalar conservation laws, although our numerical experiments in \Cref{sec:numerical} include one-dimensional systems. 

From total temporal period of observations $[0,T]$, we specify a training period
\begin{equation}
    \label{eq:trainingperiod}
    \mathcal{D}_{train}= [0,T_{L_{train}}],
\end{equation} 
where $T_{L_{train}} = L_{train}\Delta t$ and $0 < L_{train} < J$. Consistent with \eqref{model: observation} for $H = I$ (identity), we define the observable information as
\begin{equation}
\label{eq: trajectory data}
\by^v(t_j) \in {\mathbb R}^{n}, \ j=0, \cdots, L_{train},
\end{equation}
where $n$ denotes the number of spatial grid points. 
To enable structure-preserving learning, we apply the inverse transformation $\Phi^{-1}$ to obtain the corresponding sequence of conserved variable snapshots, which serve as input for surrogate training. That is, we compute
\begin{equation}
\label{eq: trajectory data - conserved}
\by^u(t_j) = \Phi^{-1}\left(\by^v(t_j)\right) \in {\mathbb R}^{n}, \ j=0, \cdots, L_{train}.
\end{equation}
We note that unlike the given observations  in \eqref{eq: trajectory data}, the noise structure in \eqref{eq: trajectory data - conserved} is neither additive nor Gaussian.  Despite this complication, the ESCFN surrogate is designed to be robust to such noise distortions and remains capable of learning stable and physically consistent dynamics from the transformed trajectory.

We emphasize that the goal of the surrogate model is to learn an approximation of the discrete-time flow map $\bPsi$ based only on the pre-processed observed trajectory $\{ \by^u(t_j)\}_{j=0}^{L_{train}}$. This learned model will then be used in the prediction step of a structurally informed ensemble Kalman filter to reconstruct and predict the evolution of $v$ in the absence of known governing equations. In so doing, we are able to leverage physical measurements while still working entirely within the space of conserved quantities in the training process -- a space that is naturally compatible with the numerical structure of hyperbolic conservation laws. The next subsections will describe how we construct the surrogate model.



\subsection{Surrogate model training}
\label{subsec: surrogate model training}

We adopt the ESCFN framework introduced in \cite{Liu2024} to build a data-driven surrogate that respects the conservative structure of hyperbolic PDEs. The scaffolding is the entropy-stable, second-order, and non-oscillatory Kurganov-Tadmor (KT) scheme \cite{KTscheme2000} applied in semi-discrete conservation form. Since the true flux function in the system is generally not known apriori, a neural network is used to approximate it in a structure-preserving way. 

The KT scheme update for the semi-discrete conservation law is given by
\begin{equation}
\frac{d}{d t} u_i(t)=-\frac{\cH_{i+1/2}(t)-\cH_{i-1 / 2}(t)}{\Delta x},
\label{eq: Kurganov Tadmor Scheme}
\end{equation}
where $\cH_{i+1/2}$ is the numerical flux across the interface between cells $i$ and $i+1$. In the ESCFN method, these fluxes are modeled by a neural network-informed expression
\begin{equation}
\begin{aligned}
\cH^{NN}_{i+1 / 2}(t):=&\frac{\cF^{\btheta}\left(u_{i+1 / 2}^{+}(t)\right)+\cF^{\btheta}\left(u_{i+1 / 2}^{-}(t)\right)}{2}
-\frac{a^{NN}_{i+1 / 2}(t)}{2}\left[u_{i+1 / 2}^{+}(t)-u_{i+1 / 2}^{-}(t)\right],
\end{aligned}
\label{eq: Numerical flux NN}
\end{equation}
where $\cF^\btheta$ is a feed-forward neural network operator parameterized by $\btheta$, and $u_{i+1 / 2}^{\pm}$ are defined as (\cite{KTscheme2000})
\begin{equation}
{\bf u}_{j+1 / 2}^{+}(t)={\bf u}_{j+1}(t)-\frac{\Delta x}{2}\left({\bf u}_x\right)_{j+1}(t), \quad {\bf u}_{j+1 / 2}^{-}(t)={\bf u}_j(t)+\frac{\Delta x}{2}\left({\bf u}_x\right)_j(t),
\label{eq: upm}
\end{equation}
where the spatial derivatives $({\bf u}_x)_j$ are computed using a Total Variation Diminishing (TVD) limiter. In particular, for each component $(u_x^i)_j, 1\leq i \leq p$, of $({\bf u}_x)_j$ we have
\begin{equation*}
( {u}_{x}^i )_{j} = \psi({r}) \left( {u}_{j + 1}^i -  {u}_{j}^i \right), \quad {r} = \dfrac{{u}_{j}^i - {u}_{j-1}^i}{{u}_{j+1}^i - {u}_{j}^i},
\end{equation*}
with the minmod-type limiter $\psi(r) = \max\big( 0, \min( r, (1+r)/2, 1) \big).$
The term $a^{NN}_{i+1 / 2}$ approximates the local maximum wave speed and is defined as
\begin{equation}
a^{NN}_{i+1 / 2}(t):=\max \left\{\rho\left(\frac{\partial \cF^{\btheta}}{\partial u}\left(u_{i+1 / 2}^{+}(t)\right)\right), \rho\left(\frac{\partial \cF^{\btheta}}{\partial u}\left(u_{i+1 / 2}^{-}(t)\right)\right)\right\}, 
\label{eq: Maximum wave speed NN}
\end{equation}
where $\rho(\cdot)$ denotes the spectral radius. This formulation ensures the surrogate flux inherits properties such as numerical stability and dissipation control.

The network $\cF^\btheta$ is implemented as a fully connected feed-forward neural network operator having $M$ layers for which the input is $u_{i+1 / 2}^{\pm} \in \bbR$ in the scalar case. We characterize $\cF^\btheta$ by trainable variables collectively denoted by $\btheta$. 
Specifically,  the network parameters $\btheta$ include the weight matrices $\mathbf{W}^{(m)} \in \mathbb{R}^{d_m \times d_{m-1}}$ for each layer $m = 1, \dots, M$, where \(d_{0} = 1\) and \(d_{M}\) is the number of unknowns, as well as the bias vectors $\mathbf{b}^{(m)} \in \mathbb{R}^{d_m}$ for each layer $m = 1, \dots, M-1$. For each $m=1, \cdots, M-1$, we then compute
\begin{equation}\label{eq:activation}
\mathbf{h}^{(m)} = \sigma\left(\mathbf{W}^{(m)} \mathbf{h}^{(m-1)} + \mathbf{b}^{(m)}\right), 
\end{equation}
where $\mathbf{h}^{(0)} = \mathbf{v}$ and $\sigma(\cdot)$ is the activation function to be specified. Finally, the output $\mathbf{f} \in \bbR^{d_M}$ is obtained by
\[
\mathbf{f} = \mathbf{W}^{(M)} \mathbf{h}^{(M-1)}.
\]
With the approximate numerical flux $\cH^{NN}$ in \eqref{eq: Kurganov Tadmor Scheme}, the problem is then solved using the total variation diminishing third-order Runge-Kutta (TVDRK3) time integration scheme \cite{Liu94}.

Given $\by^u_j = \by^u(t_j) \in \mathbb{R}^{n}$, $j = 0,\dots, L_{train}$ in \eqref{eq: trajectory data - conserved}, and annotating as \(\mathcal{N}\) the procedure of solving \eqref{eq: Kurganov Tadmor Scheme} using TVDRK3, we obtain the form for the next step predictions solved by the ESCFN method as
\begin{equation}
    \label{eq:Ndefine}
    \by^u_{j+1} = \mathcal{N}(\by^u_{j}).
\end{equation}

The network is then optimized by minimizing the \emph{recurrent loss function}
\begin{equation}
\label{eq:recurrentloss}
\mathcal{L}\left( \btheta \right) = \dfrac{1}{L_{train}} \sum_{j=0}^{L_{train}} \left\|\mathbf{u}_{NN}\left( t_{j} ; \mathbf{\btheta}\right) - \by^u_j \right\|_2^2
\end{equation}
for the trajectory dataset, where
\begin{equation}\label{eq:neuralupdate}
    \mathbf{u}_{NN}\left( t_{j}; \mathbf{\btheta} \right) = \underbrace{\mathcal{N}\circ\ldots \circ \mathcal{N}}_{j \text{ times }}\left( \by^u_0 \right).
\end{equation}

This loss encourages accurate long-term rollout by penalizing deviations between the predicted and true trajectories over the training interval. The learned surrogate $\mathcal{N}$ can then be used in ensemble forecasting and filtering, as we describe in the next section.

\subsection{Non-intrusive structure-preserving sequential data assimilation}
With the problem setup and surrogate model established, we now present the full procedure for our non-intrusive, structure-preserving sequential assimilation framework. The process consists of three primary components: (1) a preprocessing step to translate physical measurements into conserved variables; (2) a learning step that trains a structure-preserving surrogate model on the conserved variables; and (3) a data assimilation step that incorporates new observations via SIETKF.

\textbf{Preprocessing step}. The procedure begins by collecting available prior knowledge, including the physical domain constraints (e.g, the conservation law structure), the known transformation operator between conserved and physical variables $\Phi: u \to v$, and a discrete set of snapshots of the measurable physical variables \(\{\by^v(t_j)\}_{j = 0}^{J}\). Since training will be performed in the space of conserved variables $u$, we apply the inverse observation mapping to convert physical measurements to conserved variables over the training period following \eqref{eq: trajectory data - conserved}. These pre-processed snapshots $\{\by^u_j\}_{j=0}^{L_{train}}$ then serve as the training dataset for learning the underlying dynamics. Moving forward we can treat the system as if the conserved variables are directly observed. 


\textbf{Learning step}. With the training sequence of conserved variables in hand, we now train the ESCFN. As described by a series of steps in \Cref{subsec: surrogate model training}, the surrogate model mimics the conservative numerical structure of hyperbolic PDEs and is then optimized by minimizing a recurrent loss \eqref{eq:recurrentloss} over the observed trajectory, with a user-specified number of epochs.  Details for each numerical test are provided in \Cref{sec:numerical}. 


We emphasize that we are considering the case in which the available data are noisy, limited in size, and derived from only one realization of the underlying process. As a result, the ESCFN may not achieve the same accuracy as in fully supervised settings in which clean data are used. Additionally, because the initial snapshot is itself noisy, the training process must be robust to perturbations to both initial  and cumulative prediction errors. Nevertheless, the learned surrogate provides a physically consistent prior model for the assimilation stage.

\textbf{Data assimilation step}. From the training we obtain a surrogate model $\mathcal{N}$ to evolve the underlying conserved variable. Recall that  $\Phi: u \to v$, is assumed to be invertible.  We therefore can compute the unknown physical model  in \eqref{eq: EnKF predict 1 - Neural ODE} as the composition 
\begin{equation}\label{eq:Psi}
\bPsi^\btheta = \Phi \circ \mathcal{N} \circ \Phi^{-1},
\end{equation}
where again we write $\bPsi_{NN}^\btheta$ and $\bPsi_{CFN}^\btheta$ to  denote the Neural ODE and ESCFN constructions, respectively. 
The analysis step is then performed using the SETKF, as described in \Cref{subsec: SETKF}, with the learned flow map used to replace the unknown dynamics operator $\bPsi$ in \eqref{eq: EnKF predict 1}. Rather than employing the weighting matrix in \eqref{W: cov}, here we update the ensemble using localized structural information via \eqref{eq: diag gsm}, which was shown in \cite{TLda} to improve ETKF performance in systems with spatial heterogeneity or discontinuities. 

\begin{remark}
Measurements are typically acquired as physical variables, which justifies defining the observational data $\by_j = \by^v(t_j)$ in \eqref{eq: trajectory data}.  As already discussed, to enable structural preservation, we perform surrogate model training in the space of conserved variables, where the dynamics naturally conform to the form of hyperbolic conservation laws. We then use the physical variables (via transformation matrix ${\Phi}$) to compute the weighting matrix for data assimilation, since the spatial variability of a system is most meaningful and interpretable in terms of its physical quantities. This necessitates a consistent back-and-forth transformation between conserved and physical variables throughout the pipeline to ensure both structural fidelity and practical observability.    
\end{remark}

\begin{algorithm}[h!]
\caption{Non-intrusive structure-preserving sequential data assimilation (NSSDA) }
\begin{algorithmic}
\State \textbf{Input} Initial state mean $\m_{0}$ and covariance $C_0$ in \eqref{eq: init Gaussian}, observation data  $\{\by^v(t_j)\}_{j = 0}^{J}$ in \eqref{eq: trajectory data},  transformation matrix $\Phi$, ensemble size $K$, and training domain  $\mathcal{D}_{train}$ in \eqref{eq:trainingperiod}.
\State \textbf{Output} Posterior ensembles of physical variable $\{\bv_j^{(k)}\}_{j=1}^{J}$ 
\State \textbf{Preprocessing step}
\State Use the transformation matrix $\Phi$ to obtain the training dataset $\{\by^u(t_j)\}_{j=0}^{L_{train}}$ by \eqref{eq: trajectory data - conserved}.
\State \textbf{Learning step}
\State Train on $\{\by^u(t_j)\}_{j=0}^{L_{train}}$ to get the surrogate model $\cN$ as described in \Cref{subsec: surrogate model training}.
\State \textbf{Data assimilation step}
\State Initialize ensembles $\bv^{(k)}_0 \sim \cN(\m_{0},C_0))$, $k = 1,\dots,K$. 
\For{$j=1$ to $J-1$} 
\State \textbf{Prediction} Compute
\State \quad prior ensembles $\ds \wh{\bv}^{(k)}_{j+1}= \bPsi_{CFN}^\btheta(\bv^{(k)}_j)$ in \eqref{eq:Psi} 
for $k=1, \cdots, K$; 
\State \quad sample prior mean $\wh{\m}_{j+1}$ and covariance $\wh{C}_{j+1}$ according to \eqref{eq: EnKF prior 1} and \eqref{eq: EnKF prior 2};
\State \quad  structurally informed weighting matrix $W_{j+1} = \beta \wh{S}_{j+1}$ in \eqref{eq: diag gsm}.
\State \textbf{Analysis} Compute
\State  \quad posterior mean $\m_{j+1}$ as solution to the minimization problem \eqref{eq:ETRF min}, with $\by_{j+1} = \by^v(t_{j+1})$;
\State \quad posterior covariance $C_{j+1}$ by \eqref{eq:cov_update};
\State \quad posterior ensembles $\bv_{j+1}^{(k)}$ based on \eqref{eq:zeta}.
\EndFor
\end{algorithmic}
\label{alg:mETKF}
\end{algorithm}

\Cref{alg:mETKF} summarizes the entire non-intrusive structure-preserving sequential data assimilation process.
Some remarks are in order.
\begin{remark}
\label{rem:algorithm}
In some sense one could argue that  since the learning step in time domain $[0,T_{L_{train}}]$  uses  observational data in that same time frame,  the initial time for  \Cref{alg:mETKF} might instead be $t = T_{L_{train}}$. However, although obtained from data observations,  the corresponding learned solution (obtained via transformation $\bPsi_{CFN}^\btheta$ in \eqref{eq:Psi})  does not incorporate the DA analysis step.  We therefore start the process at  time $t = 0$. 
\end{remark}
\begin{remark}
\label{rem:neuralnetinalg}
\Cref{alg:mETKF} combines the structure preserving properties of the ESCFN method for learning the dynamics operator along with the SETKF data assimilation, which penalizes the incorporation of observations based on local gradient information.  Clearly one can also replace $\Psi_{CFN}^\btheta$ with  the Neural ODE operator $\Psi_{NN}^\btheta$ as well use the standard covariance matrix \eqref{W: cov} for data assimilation, (ETKF).   In our numerical experiments we consider each combination and demonstrate that the NSSDA, which incorporates structure in both the prediction and analysis steps, consistently yielding the best solutions.
\end{remark}


\section{Numerical experiments}
\label{sec:numerical}


We organize our numerical experiments into three parts. 
In the first set of experiments we consider the 1D dam break problem governed by the shallow water equations and compare different surrogate modeling strategies (Neural ODE vs.~ESCFN) and data assimilation schemes (ETKF vs.~SETKF), to assess the effects of structure preservation in both the prediction and analysis steps. For the second set of experiments we fix the modeling and assimilation methods to ESCFN and SETKF, respectively, while varying the training conditions such as noise levels and the number of training snapshots.  Note the combination of ESCFN modeling and SETKF data assimilation makes up the NSSDA in \Cref{alg:mETKF}. To evaluate performance in a more complex setting involving shocks and fine-scale features, the last set of experiments  tests  the NSSDA   on the 1D Euler equations using the Shu–Osher problem. Our results are analyzed using the relative \(\ell^2\) error given by
\begin{equation}
    \label{eq:relative l2 error}
    \mathcal{R}(\bu(t_j), \bu_{true}(t_j)) := \frac{\left\| \bu(t_j) - \bu_{true}(t_j)\right\|_2}{\left\| \bu_{true}(t_j) \right\|_2}.
\end{equation}
Here  $\bu_{true}(t)$ is the true solution at time \(t\) generated from algorithms in PyClaw \cite{clawpack,pyclaw} with the same mesh size as observations, and we compute the discrete $\ds \ell^2$ norm as \(\left\| \bu \right\|_2 = \sqrt{\Delta x\sum_{i=0}^{n} \bar{u}_i^2}\), where $\bar{u}_i$ the cell average defined in \eqref{eq: cell_average}.

The neural flux operator $\mathcal{N}$ in  \eqref{eq:neuralupdate} is constructed using a fully connected neural network, with the input and output dimensions equivalent to the number of state variables in each system. Following the architecture proposed in \cite{chencfn}, our experiments utilize five hidden layers with 64 hidden neurons at each layer.  We use the silu activation function \(\text{silu}(x) = \text{Sigmoid}(x)\cdot\text{ReLU}(x)\) in  \eqref{eq:activation}, which is smooth and differentiable, making it suitable for computing Jacobian required in \eqref{eq: Maximum wave speed NN}.  
The network used to approximate the spectral radius in  \eqref{eq: Maximum wave speed NN} is a 2-hidden-layer, 64-hidden-neurons-per-layer neural network with activation function \(\text{ReLU}(x)\) (since second-order differentiability  is not required in this case). The network to model neural ode operator \eqref{eq: neural ODE} contains four layer with 64 neurons each layer and input/output dimension to be \(2n\) (the first and latter half correspond to \(h\) and \(hu\), respectively). The number of hidden layers and corresponding neurons in each case is heuristically chosen and guided by results in \cite{Liu2024} and \cite{neuralode} respectively.

We employ the commonly used Adam stochastic optimization method \cite{Adam} to update the neural networks' weights and biases. We train our numerical models for 100 epochs in the shallow water equation examples and 500 epochs in the Euler's equation test, with a fixed learning rate of $1\times10^{-3}$ for both tests. While the number of training snapshots, \(L_{train}\), varies by type of experiment, we fix the number of  observations $J$ at 200 for the shallow water case and 800 for the Euler's case.   All implementations are based on JAX \cite{jax}, a high-performance numerical computing library that provides automatic differentiation and GPU/TPU acceleration for Python. Furthermore, to ensure robustness, we note that none of these parameters were fine-tuned.\footnote{The complete code  is  available upon request for reproducibility purposes.}  Finally,  we use $K=100$ ensemble members in the data assimilation step, initialized from either \eqref{eq:initialeasy} for the dam break problem  or \eqref{eq:initial} for the Shu-Osher problem, each perturbed by additive i.i.d.~Gaussian noise with mean $0$ and standard deviation $0.1$.


A couple of remarks are in order.  
\begin{remark}
Here we focus on the dense observation environment, with data acquired at every spatial location. We note, however, that the surrogate model is designed to operate on a full spatial grid, and is independent of the observation environment. The mismatch for the sparse observation case  introduces additional challenges. Namely, training the surrogate model would require inputting missing data or introducing correction mechanisms to compensate for unobserved regions. Investigating how to extend the structure-preserving framework in sparsely observed environments is an important direction for future research.    
\end{remark}

\begin{remark}
\label{rem:ensemblesize}
In practical applications, the choice of ensemble size depends on dimension reduction and computational feasibility. Our primary goal is to compare the impact of different surrogate models (Neural ODE vs.~ESCFN) in the prediction step and different assimilation strategies (ETKF vs.~SETKF) in the analysis step. We therefore fix $K = 100$ throughout to isolate any potential discrepancies due to methodological differences as well as to avoid additional fine-tuning. 
Supplementary numerical experiments (not reported here) confirmed that the qualitative behavior of the results did not fundamentally change for varying values of $K$.
\end{remark}

\subsection{1D dam break problem}
\label{sub:swe}
 We consider the 1D dam break problem on a flat bottom governed by shallow water equations
\begin{equation}
\label{eq:SWE}
\begin{aligned}
h_t+(hu)_x=0, \\
(hu)_t+(hu^2+\frac{1}{2}gh^2)_x=0. 
\end{aligned}
\end{equation}
Here $h=h(x,t)$ is the depth of the water, $u=u(x,t)$ is the (depth-averaged) fluid velocity, and $g = 1$ is the acceleration constant due to gravity. The fluid is initially at rest on both sides of an infinitely thin dam located at $x=0$, and the dam break is simulated by sudden removal of the dam wall at time $t = 0$ yielding initial conditions
\begin{gather}
\label{eq:initialeasy}
h(x,0)=
\begin{cases}
h_l, \quad x<x_0 \\
h_r, \quad x>x_0
\end{cases} \qan
u(x,0)=
\begin{cases}
 u_l, \quad x<x_0 \\
 u_r, \quad x>x_0
\end{cases}.
\end{gather}
In our simulations $h_l = 3.5691196$, $h_r = 1.178673$, $u_l =  -.064667$,  $u_r =  -.045197$, and \(x_0 = .003832\).
The computational domain is fixed as $[-5,5]$, and Dirichlet boundary conditions are applied with fixed left/right states matching \eqref{eq:initialeasy}. The resulting solution consists of a shock front propagating to the right and a rarefaction wave propagating to the left \cite{LeVeque92, LeVeque02}. 

The observational data $\{\by^v(t_j)\}_{j = 0}^{J}$ in \eqref{eq: trajectory data} are generated using PyClaw \cite{clawpack,pyclaw}  (HLLE Riemann Solver). The spatial resolution is set to $n=512$, and the time step $\Delta t = .005$ is chosen to satisfy the CFL condition. The total simulation time is \(T = 200\Delta t = 1\), and we use the first $L_{train} = 20$ time steps (i.e.~$\mathcal{D}_{train} = [0,.1]$) for surrogate model training. 
Finally, the  transformation operator \(\Phi\) from physical to conserved variable is given by
\(\Phi(h, u) = [h,hu]\) and \(\Phi^{-1}(h, \nu) = [h, \nu/h] \), where \(\nu=hu\) is the momentum. This allows us to treat physical observations as input and convert them to conserved variables for model training and assimilation.

\subsubsection{The effects of surrogate models and structural information}
\label{sub:swenumerics}
We first investigate the effects of employing different choices of surrogate model $\bPsi^\btheta$ in \eqref{eq: EnKF predict 1 - Neural ODE}, the Neural ODE $\bPsi^\btheta_{NN}$ and the ESCFN $\bPsi_{CFN}^\btheta$,  when used in data assimilation.  Here, too, there are two choices: the standard ETKF approach which uses the weighting matrix $W$ defined by \eqref{W: cov}, and the SETKF, which incorporates more  structural information via $W$ in \eqref{eq: diag gsm}.  To ensure a consistent and fair comparison, we fix the observation data noise to be i.i.d. Gaussian with mean \(0\) and variance \(0.1\),
and set the number of training snapshots to \(L_{train}=10\), corresponding to a training period \(t \in \left[ 0, .05 \right]\). Both the Neural ODE model and the ESCFN model are trained by the recurrent loss function.  

\begin{figure}[h!]
\begin{center}
\includegraphics[width=0.8\textwidth]{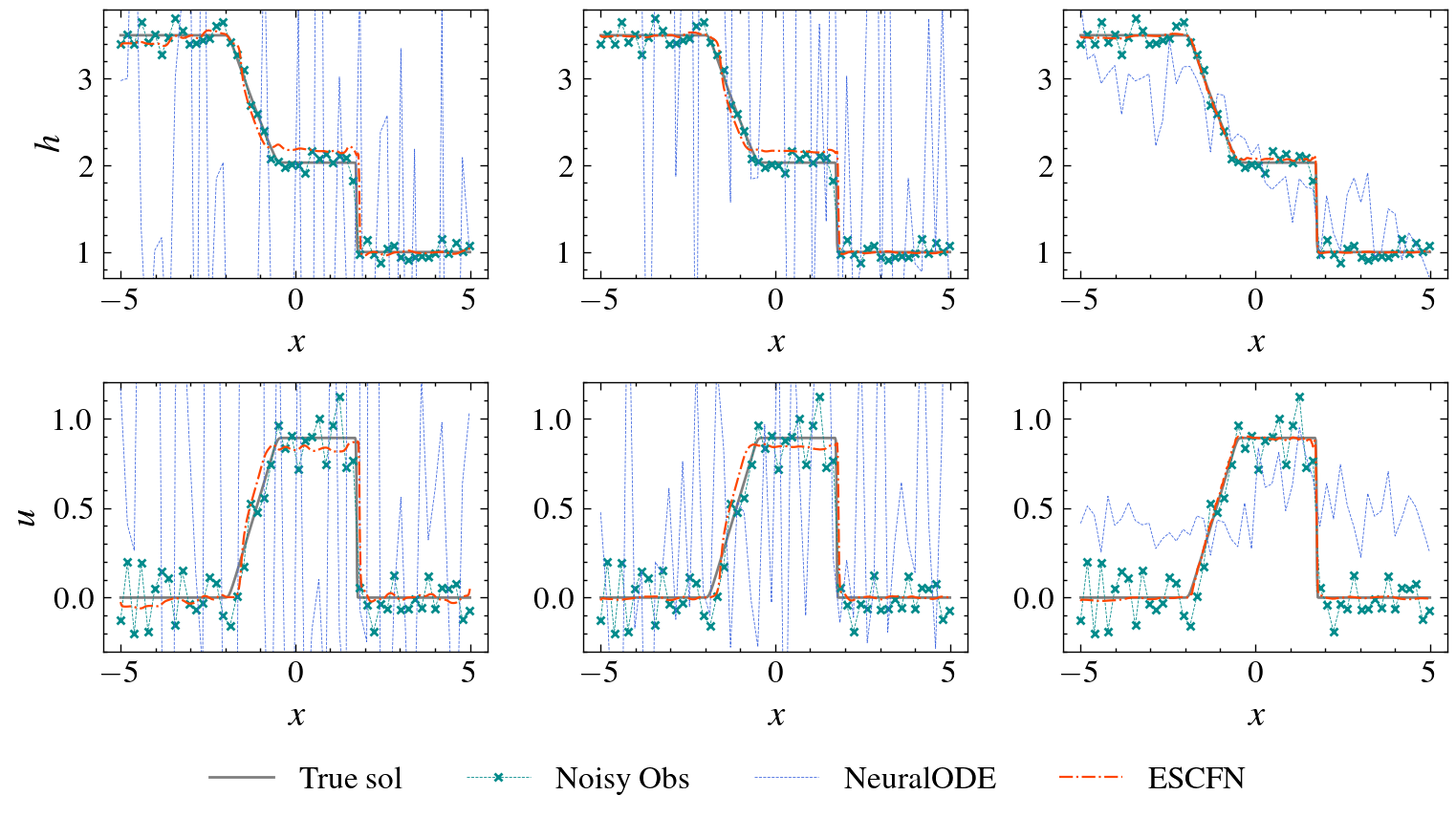}
\end{center}
\caption{The dam break problem solution at time \(t = 1\). (left) Predictions using $\bPsi_{NN}^\btheta$ and $\bPsi_{CFN}^\btheta$; (middle) the posterior ensemble means using ETKF; (right)  posterior ensemble means using SETKF. (top) height $h$; (bottom) velocity $u$.}
\label{fig: node vs ktcfn prediction}
\end{figure}

\Cref{fig: node vs ktcfn prediction} compares the model predictions using $\bPsi_{NN}^\btheta$ and $\bPsi_{CFN}^\btheta$ (left column) and their corrections via standard ETKF (middle column), and SETKF (right column), at time \(t = 1\) for the dam break problem. The top and bottom rows show the height \(h\) and velocity \(u\), respectively. It is evident that the  limited training data (i.e., \(L_{train} = 10\)) results in the Neural ODE failing to learn meaningful dynamics.  It produces highly oscillatory and inaccurate predictions (left column, blue dotted line). Moreover, applying the standard ETKF as a post-processing step does not sufficiently correct this model error (middle column). Incorporating structural information through SETKF  (right column) appears to recover structure in $h$, in spite of the highly oscillatory behavior in the prediction. 

In contrast, employing $\bPsi_{CFN}^\btheta$ in the prediction stage (red dashed line) yields essential solution properties, such as conservation and discontinuity propagation, even from limited training data. While not perfectly matching the true solution, it captures the correct qualitative trends in shock and rarefaction wave evolution, as shown in the left column. When combined with the structural information via SETKF, the result is remarkably accurate recovery of both heights \(h\) and velocity \(u\) (right column), demonstrating the synergy between structure-preserving surrogate modeling and structure-informed data assimilation.

\begin{figure}[h!]
\begin{center}
\includegraphics[width=0.6\textwidth]{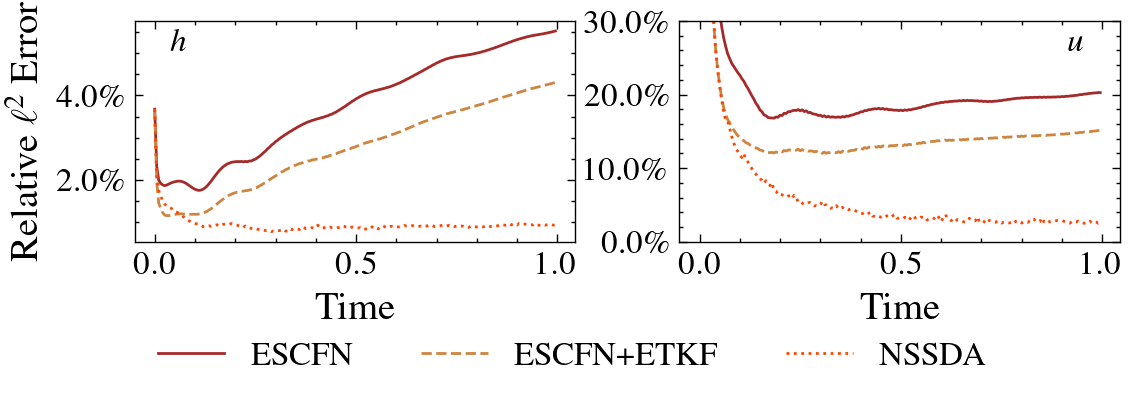}
\end{center}
\caption{Comparison of relative \(\ell^2\) errors for height \(h\) (left) and velocity \(u\) (right). The solid line shows errors for predictions using $\bPsi_{CFN}^\btheta$. The dashed and dotted lines show the errors from ETKF and SETKF post-processing, respectively.}
\label{fig: node vs ktcfn relative l2}
\end{figure}

\Cref{fig: node vs ktcfn relative l2} further quantifies the benefits of incorporating structural information into the DA analysis step via weighting matrix \eqref{eq: diag gsm} (SETKF). Specifically, it shows the relative \(\ell^2\) error of height \(h\) and velocity \(u\) at time \(t \in \left[ 0,1 \right]\), obtained from ESCFN prediction (solid line), ETKF post-processing (dashed line), and SETKF post-processing obtained as in \Cref{alg:mETKF} (dotted line). We omit the solutions obtained using $\bPsi_{NN}^\btheta$.  As already seen in \Cref{fig: node vs ktcfn prediction}, the oscillations in the Neural ODE  predictions cannot be overcome by incorporating observations. 

As expected, the  ESCFN prediction error increases over time due to cumulative model discrepancies. ETKF reduces this error significantly by incorporating observations, but the accuracy still suffers from long-term error accumulation. In contrast, SETKF leverages spatially localized discrete gradient information from the prior to better align ensemble updates with physically relevant structures, especially near discontinuities. This incorporation structural information mitigates error aggregation and results in more stable and accurate long-term predictions.

In summary, this experiment demonstrates the importance of both model structure and structural information in ensemble-based data assimilation. ESCFN provides a more accurate and stable surrogate model than a standard Neural ODE, even when the amount of training data is limited. Furthermore, incorporating structural information into the analysis step via SETKF, as provided in \Cref{alg:mETKF}, significantly improves assimilation quality and long-term prediction accuracy. We now investigate how training conditions such as data noise and training duration affect method performance.


\subsubsection{Sensitivity to noise levels and training duration}

We now examine how the observation noise magnitude and training period duration affect the performance of the NSSDA framework. To maintain consistency across comparisons, we consider several test cases with additive i.i.d.~Gaussian noise (mean \(0\) and variance \(0.1\) or \(0.2\)) in the observational data and with training durations corresponding to \(L_{train}=10, 20, 40\), i.e., training on \(\mathcal{D}_{train}=\left[ 0, .05 \right]\), \(\left[ 0, .1 \right], \text{ and } \left[ 0, .2 \right]\), respectively. In all cases, we use the ESCFN model as the surrogate prior, as Neural ODE fails to produce reliable predictions under any of these settings.

\begin{figure}[h!]
\begin{center}
\includegraphics[width=0.8\textwidth]{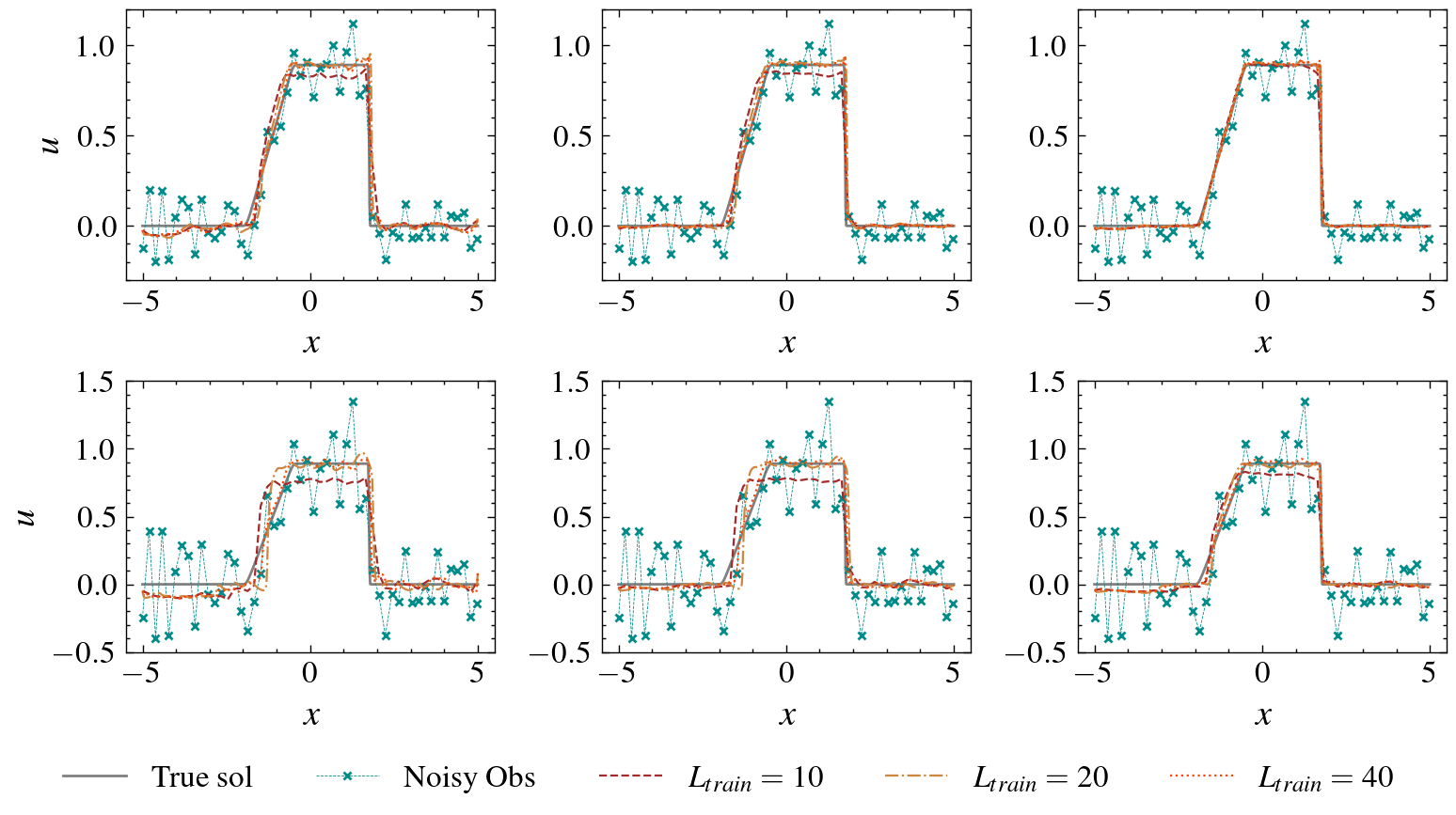}
\end{center}
\caption{(left)  ESCFN predictions with \(L_{train} = 10,20,40\); (middle) corresponding posterior ensemble means using ETKF; and (right) corresponding posterior ensemble means using SETKF at \(t = 1\). Observation variances: (top)  \(\sigma^2 = 0.1\) and (bottom) \(\sigma^2 = 0.2\).}
\label{fig: ktcfn noise trainsteps prediction u}
\end{figure}

\Cref{fig: ktcfn noise trainsteps prediction u} shows the predictions of velocity $u$ at $t=1$ for different training durations and noise levels. As expected, increasing the noise level deteriorates the ESCFN prediction quality, particularly around discontinuities. Extending the training window significantly improves model performance, however, especially in shock and rarefaction wave regions (see left column). The benefit of applying data assimilation is also clear, as  ETKF and SETKF both provide meaningful corrections to the prior predictions for $\sigma^2 = 0.1,0.2$.  SETKF consistently produces sharper reconstructions by better aligning ensemble updates with the underlying physical structure. The relative errors for height behave similarly and are omitted.

\begin{figure}[h!]
\begin{center}
\includegraphics[width=0.8\textwidth]{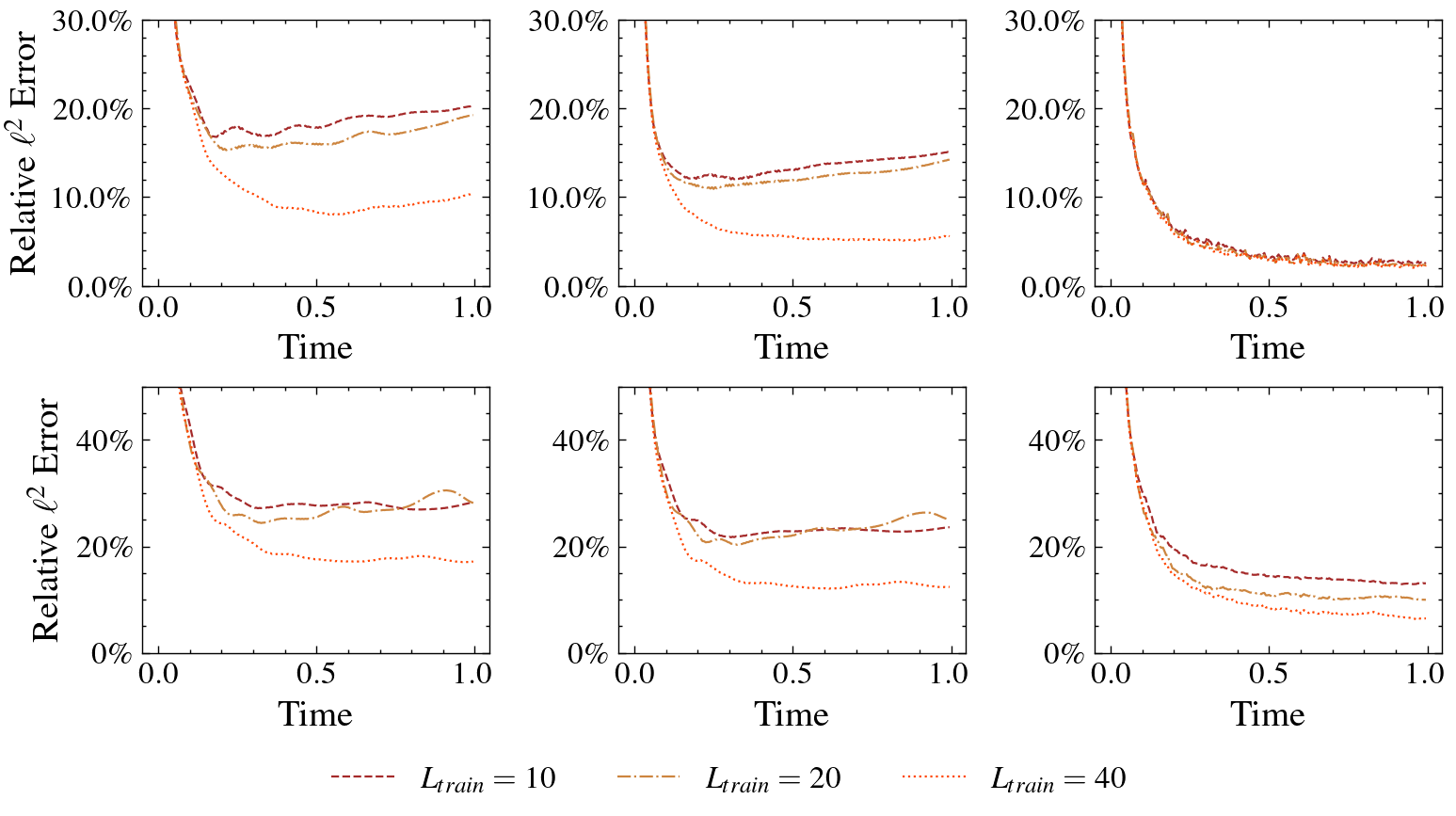}
\end{center}
\caption{(left) relative \(\ell^2\) errors for velocity \(u\) of ESCFN predictions for \(L_{train} = 10,20,40\); (middle) corresponding posterior relative $\ell^2$ error using ETKF, and (right) corresponding posterior relative $\ell^2$ error using  SETKF  over time \(t \in  \left[ 0,1 \right]\). Observation variances (top) \(\sigma^2 = 0.1\) and (bottom) \(\sigma^2 = 0.2\).}
\label{fig: ktcfn relative l2 u}
\end{figure}

\begin{figure}[h!]
\begin{center}
\includegraphics[width=0.8\textwidth]{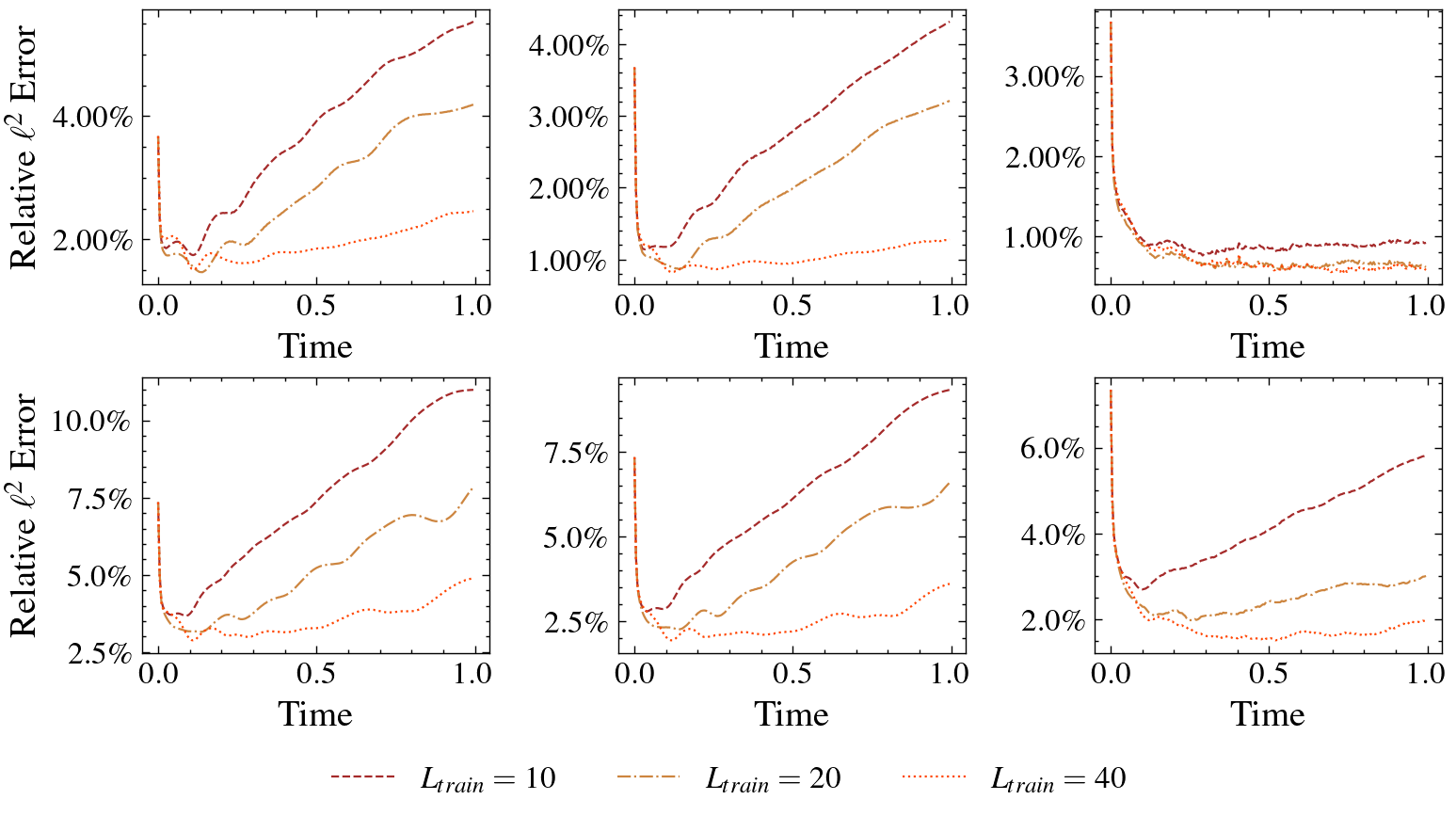}
\end{center}
\caption{(left) relative \(\ell^2\) errors for height \(h\) of ESCFN predictions for \(L_{train} = 10,20,40\); (middle) corresponding posterior relative $\ell^2$ error using ETKF, and (right) corresponding posterior relative $\ell^2$ error using  SETKF  over time \(t \in  \left[ 0,1 \right]\). Observation variances (top) \(\sigma^2 = 0.1\) and (bottom) \(\sigma^2 = 0.2\).}
\label{fig: ktcfn relative l2 h}
\end{figure}

\Cref{fig: ktcfn relative l2 u} and \Cref{fig: ktcfn relative l2 h} respectively display the $\ell^2$ errors for velocity $u$ and height $h$ for the same set of parameters over  time domain $[0,1]$.  The same trends are evident.  Noisier observations require longer training durations to control error aggregation.
Put another way, we observe that the ESCFN prior becomes increasingly informative with longer training. 
Moreover, when training data is less noisy ($\sigma^2 = 0.1$), we see that  increasing \(L_{train}\) improves accuracy for height $h$, but not for velocity $u$. 
In all we see that our NSSDA framework is robust, especially in its ability to compensate for imperfect priors under moderate levels of observation noise.

These results highlight the interplay between data quality, training duration, and structural modeling in achieving accurate long-term state estimates. When observational noise is moderate, our NSSDA framework performs robustly even with limited training. As noise increases, having access to longer temporal trajectories becomes essential for the surrogate to capture key dynamics. Importantly, the incorporation of structural information in both the surrogate model and the assimilation step enhances stability and reduces the reliance on extensive training data, making the method particularly valuable in resource-constrained or high-noise settings.


\subsection{1D Shu-Osher problem}
\label{sub:shuosher}
We now consider the Shu-Osher problem, a canonical benchmark for the one-dimensional compressible Euler equations. This problem is widely used to evaluate numerical schemes under conditions that feature both shock propagation and small-scale oscillations. The system of 1D Euler equations in conservation form is expressed by
\begin{equation}
    \begin{aligned}
        \rho_{t} + \left( \rho u \right)_{x} &= 0, \\
        \left( \rho u \right)_{t} + \left( \rho u^2 + p \right)_{x} &= 0, \\
        \left( E \right)_{t} + \left( u (E + p) \right)_{x} &= 0,
    \end{aligned}
    \label{eq: Euler equation}
\end{equation}
where $\rho(x,t)$ is the fluid density, $u(x,t)$ is the velocity, $p(x,t)$ is the pressure, and $E(x,t)$ is the total energy, defined by $E=\frac{p}{\gamma-1} + \frac{1}{2} \rho u^2$, with $\gamma=1.4$ being the adiabatic constant.

We consider the spatial domain $(-5,5)$ with Dirichlet boundary conditions and initial conditions given by
\[
\begin{aligned} \label{eq:initial}
& \rho(x, 0)= \begin{cases}3.857135, & \text { if } x \leq -4, \\
1+0.2 \sin (5 x), & \text { if } -4<x \leq x_1, \quad u(x, 0)= \begin{cases}2.62936, & \text { if } x \leq -4, \\
0, & \text { otherwise, }\end{cases} \\
1+0.2 \sin (5 x) e^{-\left(x-x_1\right)^4} & \text { otherwise, }\end{cases} \\
& p(x, 0)=\left\{
\begin{array}{ll}
 10.33333, & \text { if } x \leq -4, \\
1., & \text { otherwise, }
\end{array} \quad E(x, 0)=\frac{p_0}{\gamma-1}+\frac{1}{2} \rho(x, 0) u(x, 0)^2 . \right.
\end{aligned}
\]
with $x_1=3.29867$ and $\gamma=1.4$. The trajectory is generated by employing PyClaw (HLLE Riemann Solver) with  \(\Delta t = .002\) to solve the system. The observations are then obtained via \eqref{model: observation}  with $H = I$  and variance $\sigma^2 = 0.2$. 
The training time window $\mathcal{D}_{train}$ in  \eqref{eq:trainingperiod} is taken as $[0,.04]$ (i.e., \(L_{train} = 20\)), while the trajectory for data assimilation is evolved up to \(T = 1.6\) (i.e., \(J = 800\) steps). The spatial grid size is chosen so that \(n = 512\). Finally, we note that the transformation operator \(\Phi\) from physical to conserved variables is defined by 
\(\Phi(\rho, u, p) = \left[\rho,\rho u, \frac{p}{\gamma-1}+\frac{1}{2} \rho u^2\right]\), with inverse mapping \(\Phi^{-1}(\rho, \nu, E) = \left[\rho, \frac{\nu}{\rho},  \left( \gamma - 1 \right)\left( E - \frac{\nu^2}{2\rho} \right) \right] \), where $\nu = \rho u$ denotes the momentum.


\begin{figure}[h!]
\begin{center}
\includegraphics[width=0.7\textwidth]{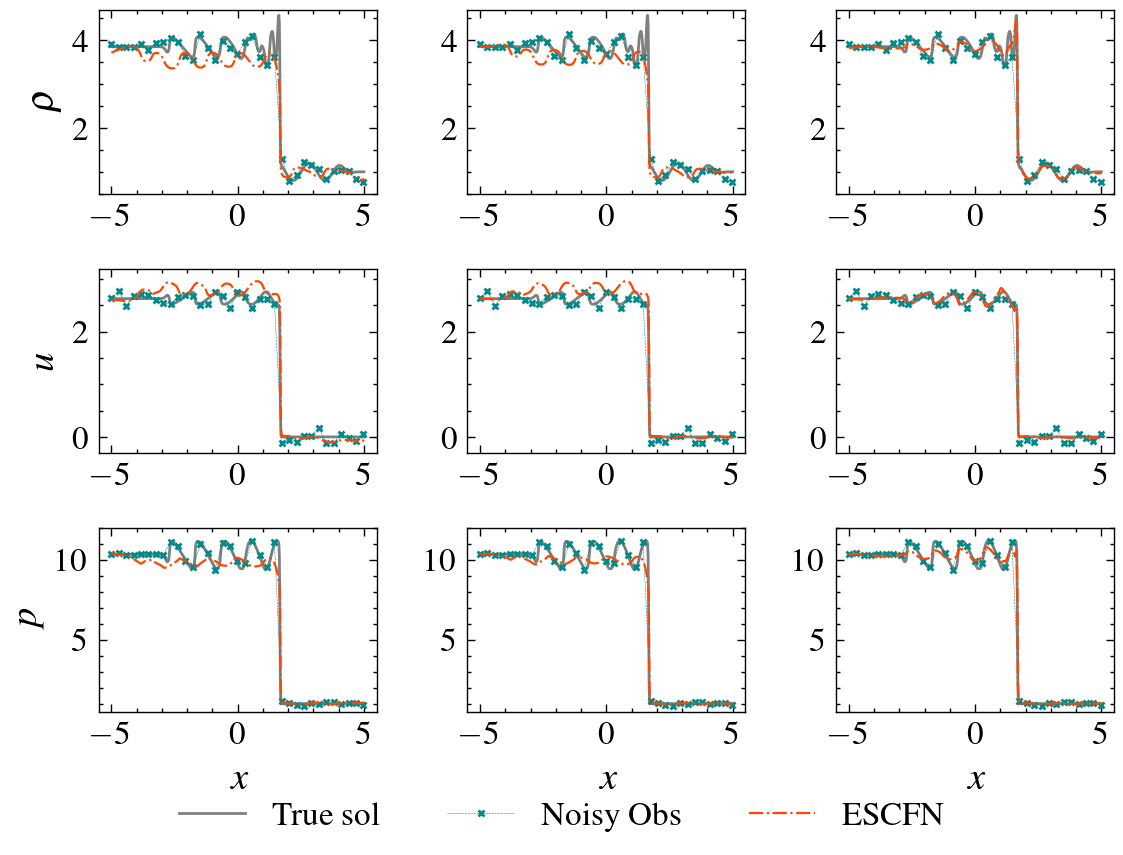}
\end{center}
\caption{(left) ESCFN  predictions; (middle column) corresponding posterior ensemble means using ETKF; and (right) corresponding posterior ensemble means using SETKF  at \(t = 1.6\). (top) density $\rho$; (middle row) velocity $u$;  and (bottom)  pressure \(p\).}
\label{fig: euler ktcfn predictions}
\end{figure}

Similar to \Cref{fig: ktcfn noise trainsteps prediction u} for the 1D dam break problem, \Cref{fig: euler ktcfn predictions} compares the qualitative  predictions using ESCFN to learn the dynamics of height \(h\), velocity \(u\) and pressure \(p\) at time $t=1.6$  (left), the corresponding correction using ETKF (middle column) and the corresponding correction using SETKF, or equivalently, \Cref{alg:mETKF} (right).  The ESCFN model prediction captures the main wave structures but underestimates density $\rho$ and pressure $p$, while overestimating velocity $u$ in oscillatory regions. Post-processing with ETKF yields only modest improvements over the prior, with slight corrections observed in the oscillatory regions of density and in the smoother right tail of the velocity field. In contrast, NSSDA recovers significantly finer details across all variables, particularly in regions with complex wave interactions, by effectively aligning ensemble updates with local structural features extracted from the prior.

\begin{figure}[h!]
\begin{center}
\includegraphics[width=0.8\textwidth]{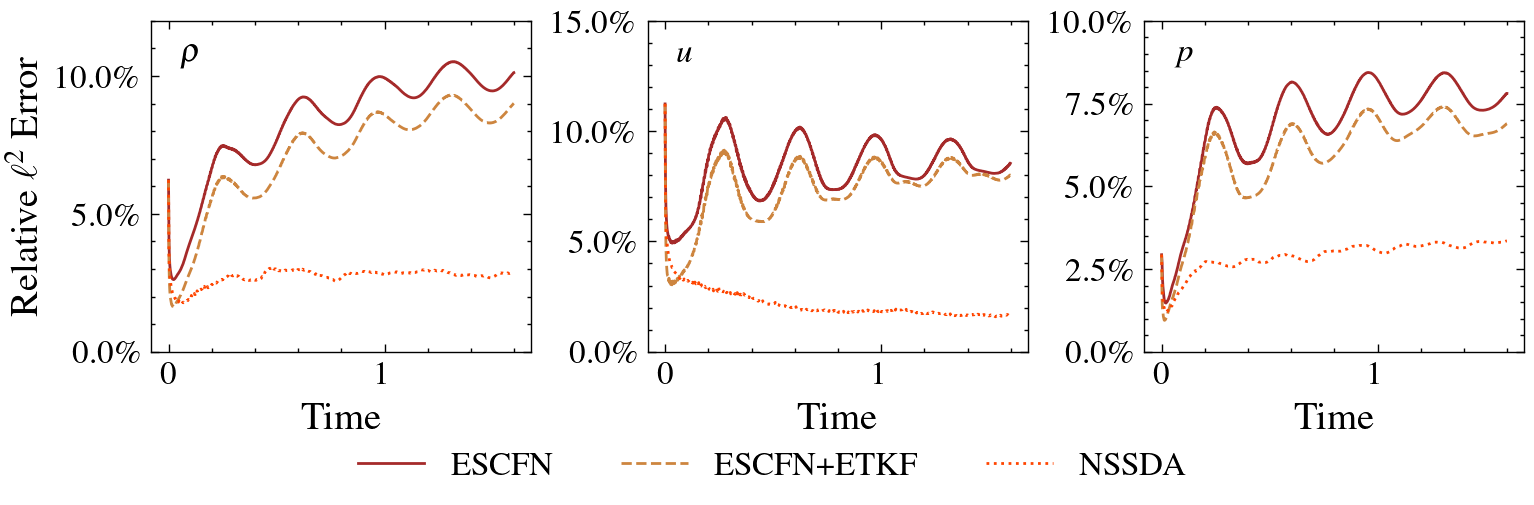}
\end{center}
\caption{(left) relative \(\ell^2\) errors for density \(\rho\); (middle) velocity \(u\); and (right) pressure \(p\) in the Shu-Osher test case. The solid, dashed, and dotted line shows errors for ESCFN predictions, ETKF and SETKF post-processing, respectively.}
\label{fig: euler ktcfn relative l2}
\end{figure}

These observations are further validated in \Cref{fig: euler ktcfn relative l2}, which reports the relative $\ell^2$ errors over time. Across all state variables, NSSDA consistently outperforms both the ESCFN prior and the standard ETKF correction. Similar to the dam break problem case, here we see that the ESCFN prior  suffers from growing errors due to accumulated model inaccuracies. While ETKF provides only limited improvement, SETKF successfully stabilizes the error growth. This highlights the advantage of incorporating localized structural information into the assimilation process, enabling more accurate and reliable recovery in complex dynamical regimes.


\section{Concluding remarks}
\label{sec:conclusion}
Our novel  non-intrusive structure-preserving sequential data assimilation (NSSDA) framework integrates a physics-informed surrogate model via the entropy stable conservative flux form neural net (ESCFN), with the  structurally informed ensemble transform Kalman filter (SETKF). This framework is designed to address the challenge of model-free forecasting and correction in hyperbolic conservation laws using only observed data without explicit knowledge of the underlying PDEs.

Our numerical experiments demonstrate that structural information plays a critical role in both the surrogate modeling and assimilation steps. Compared to Neural ODEs, the ESCFN surrogate successfully captures key features such as shocks and rarefaction waves, even under limited and noisy training data. Moreover, incorporating gradient-based structural information into the ensemble update step via SETKF yields substantial improvements in long-term prediction accuracy and stability, especially in regions exhibiting oscillatory or discontinuous behavior. Moreover, the NSSDA framework is robust to observational noise and can effectively recover accurate states across varying training durations. In particular, SETKF provides consistent gains by reducing error growth, stabilizing ensemble updates, and aligning corrections with physically relevant structures.

This work can be extended in multiple ways. A natural next step is to extend the NSSDA framework to higher-dimensional systems, particularly those involving complex geometries and boundary conditions. Another important challenge lies in adapting the framework to settings with sparse or indirect observations, where data input or correction strategies may be necessary to compensate for unobserved regions. Additionally, this work assumes a known transformation operator between conserved and physical variables. This operator is often unknown or only partially understood in practical applications, however, and  raises the important research question regarding whether and how underlying conservation structures can be inferred directly from observed data. Incorporating domain knowledge into the learning of such transformations may significantly improve the framework’s flexibility and generalizability. These extensions are essential steps toward applying structure-informed, model-free data assimilation to complex real-world systems in geophysics, climate science, and beyond.

\bibliography{mybibfile}

\end{document}